\documentclass[12pt,a4paper]{article}
\usepackage{enumerate}
\usepackage{amsfonts,amsmath}
\usepackage{graphicx,amssymb}
\usepackage[bookmarks=true]{hyperref}
\newcommand {\IZ}{\mathbb{Z}}
\newcommand {\IN}{\mathbb{N}}  %N
\newcommand {\IR}{\mathbb{R}}   %R

\newcommand {\N}{\mbox{${\cal N}$}}

\newcommand {\w}{\mbox{${\omega}$}}

\newcommand{\proof}{\noindent {\sc Proof.\ \ }} %PROOF

\newtheorem{stat}{Statement}
\newtheorem{cor}[stat]{Corollary}
\newtheorem{thm}[stat]{Theorem}
\newtheorem{lemma}[stat]{Lemma}
\newtheorem{remark}[stat]{Remark}

\newcommand{\bd}{\begin{displaymath}}
\newcommand{\be}{\begin{equation}}
\newcommand{\beq}{\begin{eqnarray}}
\newcommand{\ba}{\begin{array}}
\newcommand{\ed}{\end{displaymath}}
\newcommand{\ee}{\end{equation}}
\newcommand{\eeq}{\end{eqnarray}}
\newcommand{\ea}{\end{array}}

\newcommand{\indicator}{\mbox{\large\bf$1$}}

\newcommand{\abs}[1]{\left|#1\right|}
\newcommand{\Prob}{{\rm I\hspace{-0.8mm}P}}

\newcommand{\dt}{\partial_t}
\newcommand{\dx}{\partial_x}
\begin{document}
\title{Macroscopic stability for nonfinite range kernels}
\author{
T.S. Mountford$^{1,4}$,
K.\ Ravishankar$^{2,4}$,
E.\ Saada$^{3,4}$,
}
\date{September 16, 2009}
\maketitle
$$ \ba{l}
^1\,\mbox{\small
Institut de Math\'ematiques, \'Ecole Polytechnique F\'ed\'erale,
}\\
$\quad$ \mbox{\small
Station 8, 1015 Lausanne, Switzerland.
{\sl  thomas.mountford@epfl.ch}
}\\
^2\, \mbox{\small Department of Mathematics, SUNY,
College at New Paltz,
}\\
$\quad$ \mbox{\small
NY, 12561, USA.
{\sl  ravishak@newpaltz.edu}
}\\
^3\, \mbox{\small CNRS,  UMR 6085, LMRS, Universit\'e de Rouen,
}\\
$\quad$ \mbox{\small
 76801 Saint Etienne du Rouvray, France.
{\sl ellen.saada@univ-rouen.fr}
}\\
^4\, \mbox{\small Institut Henri Poincar\'e, Centre Emile Borel,
}\\
$\quad$ \mbox{\small
11 rue Pierre et Marie Curie, 75005 Paris, France.
}\\
\ea
$$

\begin{abstract}
\noindent
We extend the strong macroscopic stability introduced in
Bramson \& Mountford (2002) for one-dimensional asymmetric
exclusion processes with finite range to a large class of one-dimensional
conservative attractive models (including misanthrope
process) for which we relax the requirement of
finite range kernels. A key motivation is the extension
of  constructive hydrodynamics result of
Bahadoran {\em et al.} (2002, 2006, 2008) to
nonfinite range kernels.
\end{abstract}

\vskip 1 truein
\eject

\section{Introduction}
In this note we consider a general class of (at least potentially)
long range one dimensional conservative attractive particle  systems
(which will be shortly specified). The paper is motivated by the
recent series of papers \cite{bgrs1,bgrs2} and \cite{bgrs3}.
Here the hydrodynamic limits of various systems was established. The
needed conditions were extremely general, to the point where it was
not necessary to suppose that a full characterization of translation
invariant equilibria had been established.  Briefly the argument
built on the approach of \cite{av} which establishes hydrodynamic
limits for   Riemannian   initial profiles.  Then a general argument
was given to pass from this particular case to general initial
profiles.  A key part of this passage was the existence of a
macroscopic stability criterion for the particle systems whereby the
known behaviour of a system corresponding to a step-function profile
could yield information about systems corresponding to more general
(but close) initial profiles.

We now detail the processes involved.
The state space is ${\bf X}=\{0,\cdots,K\}^{\IZ}$. The evolution
consists in particles' jumps, according to the generator
\be \label{generator}
Lf(\eta)=\sum_{x,y\in{\IZ}}p(y-x)b(\eta(x),\eta(y)) \left[
f\left(\eta^{x,y} \right)-f(\eta) \right] \ee
for a local function $f$, where $\eta^{x,y}$ denotes the new state
after a particle has jumped from $x$ to $y$ (that is
$\eta^{x,y}(x)=\eta(x)-1,\,\eta^{x,y}(y)=\eta(y)+1,\,
\eta^{x,y}(z)=\eta(z)$ otherwise), $p$ is the particles' jump
kernel, that is $\sum_{z\in\IZ}p(z)=1$, and $b\ :\
\IZ^+\times\IZ^+\to \IR^+$ is the jump rate. We assume that $p$ and
$b$ satisfy :
\begin{tabbing}
{\em (A1)} \=
The greatest common  divisor of the set $\{x:p(x) \ne 0\}$ equals 1\\
\> (irreducibility);\\
{\em (A2)} \> $p$ has a finite first moment,
that is  $\mu_1 \ = \sum_{z\in\IZ}\abs{z} p(z)<+\infty$,\\
\>
and a positive mean,
that is  $0<\mu \ = \sum_{z\in\IZ}z p(z)$;\\
{\em (A3)} \> $b(0,.)=0,\,b(.,K)=0$ (no more than $K$
particles per site), and\\
\> $b(1,K-1)>0$;\\
{\em (A4)} \> $b$ is nondecreasing (nonincreasing) in its first
(second) argument\\
\> (attractiveness).
 \end{tabbing}
 For us the departure from the previous works mentioned is in
assumption {\em (A2)}, which replaces the ``finite range''
assumption.

Let $\mathcal I$ and $\mathcal S$ denote respectively the set of
invariant probability measures for $L$, and the set of
shift-invariant probability measures on $\mathbf{X}$. It was derived in
\cite[Proposition 3.1]{bgrs2} that
\be \label{characterize_2} \left({\mathcal I}\cap{\mathcal
S}\right)_e= \left\{\nu^\rho,\,\rho\in{\mathcal R}\right\}\ee
 with $\mathcal R$ a closed subset of $[0,K]$ containing $0$
and $K$, and $\nu^\rho$ a shift-invariant measure such that
$\nu^\rho[\eta(0)]=\rho$. (The index $e$ denotes extremal elements.)
The measures $\nu^\rho$ are stochastically ordered:
$\rho\leq\rho'\Rightarrow\nu^\rho\leq\nu^{\rho'}$.
\\

The result to be announced in the next section considers
``naturally'' coupled systems.  It is time to detail the natural
coupling in force throughout this paper. We suppose given on a space
$\Omega$ a family of independent marked Poisson processes $\N^{x,y}
$ of rate $p(y-x)||b||_\infty$ where $||b||_\infty = \max_{0 \leq
i,j \leq K}\{b(i,j)\}$ and associated to each point $t \in \N^{x,y}
$ are uniform random variables $U(x,y,t)$ on $[0,1]$ which are
independent over all  $(x, y)\in\IZ^2 $ and $t\ge 0$. We also
assume that   the Poisson, uniform random variables (mutually
independent and independent of the previous processes $\N^{x,y} $
and $U(x,y,t)$) that we will need for the proofs of this note are
defined on $\Omega$. We denote by $\Prob$ the probability
measure on $\Omega$. The initial configurations are defined on a
probability space $(\Omega_0,\Prob_0)$. Given an initial
configuration $\eta_0(\w_0) \in \{0, \cdots, K\}^\IZ $ and a
realization $\w$ of the Poisson processes
 and uniform random variables,  we
construct a process $(\eta_t: t \geq 0)
:=(\eta_t(\eta_0(\w_0),\omega): t \geq 0)$ by stipulating that the
process $\eta_.$ jumps from $\eta_{t^-} $ to $\eta_t =
\eta_{t^-}^{x,y} $ only if $t \in \N^{x,y} $ and $U(x,y,t) \leq
b(\eta_{t^-}(x),\eta_{t^-}(y))/||b||_\infty $.   We note that
through the above (Harris) graphical construction (see \cite {bgrs3}
for details), an evolution is constructed given any initial
configuration. Thus for any two configurations $\eta_0 $ and
$\xi_0$ we
have two naturally coupled processes, through basic coupling.  \\

We now discuss the macroscopic stability property which was introduced in \cite{bm}.
For this we introduce some notation.
For two bounded measures $\alpha(dx)$, $\beta(dx)$ on $\IR$
with compact support, we define
\be\label{def_delta}
\Delta(\alpha,\beta):=\sup_{x\in\IR}\abs{
\alpha((-\infty,x])-\beta((-\infty,x])}.
\ee
Let $N\in\IN$ be the scaling parameter for the
hydrodynamic limit, that is the inverse of the macroscopic distance
between two consecutive sites. Let
\[
\alpha^N(\eta)(dx)=N^{-1}\sum_{y\in\IZ}\eta(y)\delta_{y/N}(dx)
\in{\mathcal M}^+(\IR) \]
denote the empirical measure of a configuration $\eta$ viewed on
scale $N$, and ${\mathcal M}^+(\IR)$   denote   the set of positive
measures on $\IR$ equipped with the metrizable topology of vague
convergence, defined by convergence on continuous test functions
with compact support.

By macroscopic stability we mean that
$\Delta$ is an ``almost'' nonincreasing
functional for a pair of coupled evolutions
$(\eta_t,\xi_t: t \geq 0)$ where $\eta_0$ and $\xi_0$
are any two configurations
 with a finite number of particles, in the following sense. There
 exist constants $C>0$ and
$c>0$, depending only on $b(.,.)$ and $p(.)$, such that  for
every $\gamma>0$, the event
\be\label{bmevent} \forall t>0:
\,\Delta(\alpha^N(\eta_t(\eta_0,\omega)),\alpha^N(\eta_t(\xi_0,\omega)))
\leq \Delta(\alpha^N(\eta_0),\alpha^N(\xi_0))+\gamma
\ee
has $\Prob$-probability at least
$1-C(\abs{\eta_0}+\abs{\xi_0})e^{-cN\gamma}$, where
$\abs{\eta}:=\sum_{x\in\IZ}\eta(x)$.
\\
The strong macroscopic stability property was introduced in
\cite[Section 3]{bm} to determine the existence of stationary blocking measures
for one-dimensional exclusion processes with a random walk kernel $p(.)$
having finite range and positive mean.
It was then applied to  models considered in this note in \cite{bgrs1,bgrs2,bgrs3} with the additional assumption that the jumps had a finite range.  An essential ingredient for this property is the attractiveness of the model.

\begin{remark}
While in \cite{bm} and in the rest of this note a function $\Phi$ is used (see \eqref{Phi} below) to measure distance between configurations, we use $\Delta$ in the discussion above since it is more appropriate for hydrodynamics.
An elementary computation shows that the statement in \eqref{generalized3.1} remains unchanged
whether one uses  $\Phi$ or $\Delta$.
\end{remark}

In Section \ref{result} we state the macroscopic stability result, and its application to strong hydrodynamics. In Section \ref{discrepancies} we prove it, through an analysis of the evolution of labelled discrepancies. Section \ref{macrostab2} is devoted to two properties needed for hydrodynamics of the particle system.

\section{The result}\label{result}
We fix
\be\label{choice_L}
L > 10 (\mu_1 +1).
\ee
\begin{thm}\label{th:generalized_BM_3.1}
Let $\eta^i_., \ i= 1,2$ be two processes  both
generated by the  same Harris system  with initial configurations $\eta^i_0,\ i= 1,2$  such that
\be\label{1}
\sum_{|x| \geq LN} (\eta^1_0(x) + \eta^2_0(x)) \quad   = \quad 0.
\ee
We set, for $t\geq 0$,  $x \in \IZ$,
\be\label{Phi} \Phi_t(x) \ = \ \sum_{y\geq x} (\eta_{t}^1(y) - \eta_{t}^2 (y)).\ee
Then, for each $\epsilon >0$,
\be\label{generalized3.1}
\Prob\left(\sup_{x \in \IZ}\Phi_t(x)- \sup_{x \in \IZ}\Phi_0(x)> \epsilon N\right) \leq Ce^{-cN}
\ee
 for all   $t \in [0,N]$   and $N$, and appropriate $c>0$ and $C$, depending on $\epsilon $ and $ L $ but not on
$N$ or $\eta^i_0,\ i= 1,2$.
\end{thm}
 One can   extend this type of result to initial joint configurations
which agree outside interval $(-LN,LN)$ but
do not necessarily satisfy Condition \eqref{1} by an approach which relies on Theorem \ref{thm2}, Section 3.\\

 Theorem \ref{th:generalized_BM_3.1}  has practical consequences to hydrodynamics.
It enables us to extend the hydrodynamics derived in \cite{bgrs1,bgrs2,bgrs3},
for which the assumption
$p(.)$ finite range (that is there exists $M>0$ such
that $p(x)=0$ for all $|x|>M$) was needed, to any transition kernel $p(.)$ satisfying
{\em (A1)}, {\em (A2)}. We now state this hydrodynamic result in a more general form,
namely a strong hydrodynamic limit (which was the setup in \cite{bgrs3}).
\begin{thm}\label{th:hydro-extended}
Assume $p(.)$ has a finite  third moment  $\mu_3=\sum_{z\in\IZ} |z|^3p(z)<+\infty$.
 Let $(\eta^N_0,\,N\in\IN)$ be a
sequence of $\bf X$-valued random variables on $\Omega_0$. Assume
there exists a measurable $[0,K]$-valued profile $u_0(.)$ on $\IR$
such that
\be\label{initial_profile_vague}
\lim_{N\to\infty}\alpha^N(\eta^N_0)(dx)=
u_0(.)dx,\quad\Prob_0\mbox{-a.s.}\ee
that is,
\[
\lim_{N\to\infty}
\int_\IR\psi(x)\alpha^N(\eta^N_0)(dx) =\int
\psi(x)u_0(x)dx,\quad\Prob_0\mbox{-a.s.} \]
for every continuous function $\psi$ on $\IR$ with compact support.
 Let $(x,t)\mapsto u(x,t)$ denote the unique entropy
solution to the scalar conservation law
\be \label{hydrodynamics} \dt u+\dx[G(u)]=0 \ee
with initial condition $u_0$, where $G$ is a Lipschitz-continuous
flux function  (defined in \eqref{flux} below) determined by $p(.)$
and $b(.,.)$.
Then,  with $\Prob_0\otimes\Prob$-probability one, the convergence
\be \label{later_profile}
\lim_{N\to\infty}\alpha^N(\eta_{Nt}(\eta^N_0(\omega_0),\omega))(dx)=u(.,t)dx
\ee
holds uniformly on all bounded time intervals. That is, for every
continuous function $\psi$ on $\IR$ with compact support, the
convergence
\[
\lim_{N\to\infty}
\int_\IR\psi(x)\alpha^N(\eta^N_{Nt})(dx) =\int \psi(x)u(x,t)dx\]
holds uniformly on all bounded time intervals.
\end{thm}
While this condition on kernel $p(.)$ is probably nonoptimal, we have chosen not to pursue this question, prefering to give a simple and direct argument.
We recall from \cite[pp.1346--1347 and Lemma 4.1]{bgrs2} the
definition of the Lipschitz-continuous macroscopic flux function
$G$. For  $\rho\in\mathcal R$, let \be \label{flux}
G(\rho)=\nu^\rho\left[\sum_{z\in\IZ}zp(z)b(\eta(0),\eta(z))
\right]; \ee
 this represents the expectation, under the shift invariant
equilibrium measure with density $\rho$, of the microscopic current
through site $0$. On the complement of $\mathcal R$, which is at
most a countable union of disjoint open intervals, $G$ is
interpolated linearly.
 A  Lipschitz constant $V$ of $G$ is determined by the
rates  $b(.,.),p(.)$  in \eqref{generator}:
\[
V=2\mu_1\sup_{0\le a\le K,0\le k<
K}\{b(a,k)-b(a,k+1),b(k+1,a)-b(k,a)\}. \]

To obtain the above theorem by a constructive approach, one proceeds by first proving hydrodynamics for Riemann initial profiles and then by a general argument motivated by Glimm scheme obtain the general hydrodynamics by an approximation scheme. We now explain briefly how this approximation result is proved in the setup of \cite{bgrs3}, that is  $\Prob_0\otimes\Prob$-a.s.  convergence. Therefore all the involved processes are evolving on a common realization $(\w_0,\w)\in\Omega_0\times\Omega$, that we omit from the notation for simplicity. This proof is based on an interplay of macroscopic properties for the conservation law and microscopic properties for the particle system, in particular macroscopic stability and finite propagation property, both valid at microscopic as well as at macroscopic level. The useful properties of the entropy solution $u(.,t)$ to the conservation law are summarized in \cite[Proposition 4.1]{bgrs3}.

For $T\in\IR^+$, the time interval $[0,T]$ is partitioned by $\{t_1,t_2, \cdots\}$ into intervals of equal length. At the macroscopic level the general profile at the beginning of each time step $t_k$ (that is, the solution  $u(.,t_k)$ of the conservation law) is approximated by a step function $v_k(.)$; the time and space steps are chosen so that the  Riemann  solutions of different spatial steps (``waves'') do not interact during $[t_k,t_{k+1}]$. Macroscopic stability for the conservation law  implies that \[\Delta(u(.,t_{k+1})dx,v_{k}(.,t_{k+1}-t_k)dx) \leq \Delta(u(.,t_{k})dx,v_{k}(.)dx)\]
where $v_k(.,t-t_k)$ is the entropy solution of the conservation law at time $t$ with initial condition $v_k(.)$ at time $t_k$. We denote by $\xi^{N,k}$ the initial configuration at time $Nt_k$ which is a ``microscopic version'' of $v_k$, and by $\xi^{N,k}_{N(t-t_k)}$ the evolved configuration at time $Nt$. By this we mean that
\[\lim_{N\to\infty}\Delta(\alpha^N(\xi^{N,k}),
v_k(.)dx)=0;\]
for $k=0$, this follows from  an ergodic theorem for the densities (notice that the measures $\nu^\rho$
are not necessarily product); for $k\ge 1$, this follows from Riemann hydrodynamics applied to a profile with constant density.
At the microscopic level,
\[\Delta (\alpha^N(\eta^N_{Nt_{k+1}}),\alpha^N(\xi^{N,k}_{N(t_{k+1}-t_k)})) \leq \Delta (\alpha^N(\eta^N_{Nt_k}),\alpha^N(\xi^{N,k}))+\epsilon\]
with probability greater than $1 - CN e^{-cN\epsilon}$ by
macroscopic stability at the particle level (that is, Theorem
\ref{th:generalized_BM_3.1}). If we know that
\[\lim_{N \to \infty} \Delta(\alpha^N(\xi^{N,k}_{N(t_{k+1}-t_k)}),v_k(.,t_{k+1}-t_k)dx)= 0\]
  then we would have shown that the error
 \[|(\Delta (\alpha^N(\eta^N_{Nt_{k+1}}),u(.,t_{k+1})dx) - \Delta (\alpha^N(\eta^N_{Nt_{k}}),u(.,t_k)dx))|\]
 is small
 and the proof can be completed by induction on $k$. This last step requires  patching together Riemann hydrodynamics
for which one needs the finite propagation property for the particle system (which requires that $p(.)$
has a finite third  moment).  The bound $CNe^{-cN \epsilon }$ is not necessary for the argument.\\

Since the ergodic theorem for densities and the finite propagation property were stated in \cite{bgrs1,bgrs2,bgrs3} for finite range transition kernels $p(.)$, we state and prove their extension
to nonfinite range kernels for the sake of completeness (see Section \ref{macrostab2}).

\section{Discrepancies}\label{discrepancies}

For two processes $(\eta^1_t : t \geq 0)$ and $(\eta^2_t : t \geq 0)$ we say that there is a discrepancy at $x \in \IZ$ at time $t$
if
$\eta^1_t(x) \ne \eta^2_t(x) $.  If
$\eta^1_t(x)- \eta^2_t(x) = h\in\IN\setminus\{0\}$
we say that there are $h$ 1/2 discrepancies at site $ x \in \IZ $ at time $t \geq 0$.  We similarly speak of 2/1 discrepancies.   Indeed, we do not permit different types of discrepancies to
share the same site.  Given condition \eqref{1}, for two processes as in  Theorem \ref{th:generalized_BM_3.1}
there are
only a finite (and, given the common Harris system, decreasing
 since the model is attractive) number of discrepancies of either type.
It will be of interest to consider the time evolution of
discrepancies; to this end we will, as in \cite{bm}, label them:
for, say, 1/2 discrepancies, we will introduce  the processes
$(X^{x,i}_t: t \geq 0)$ of their positions, for $x \in (-LN, LN)$
and $1 \leq i \leq K$, taking values in $ \IZ \cup \{ \Delta \} $
where $ \Delta $ is a graveyard site.  For 2/1 discrepancies, we
will introduce processes $(Y^y_t : t \geq 0 )$ for $y$ in some
labeling set $J$,  a cemetery state, $\Delta ^\prime$, such  that at
all times $t$, $\{z: \eta^2_t(z)> \eta^1_t(z)\} \cup \{ \Delta
^\prime \}$ is equal to the union of the positions $Y^y _t $ with
multiplicities respected,   that is
$$\sum_{z\in\IZ}(\eta^2_t(z)- \eta^1_t(z))^+\delta_{z}=\sum_{y\in J}{\bf 1}_{\{Y^y _t\not=\Delta ^\prime\}}\delta_{Y^y _t}$$
   A decrease of discrepancies corresponds to
the coalescence of a 1/2 and a 2/1 discrepancies, due to the jump of one of them to the site where the other
is; in that case, we will  make the label  of a 1/2 discrepancy
 (not necessarily the one involved in the jump, see case [e] below)  jump to $ \Delta $, and the label of the
2/1 discrepancy jump to $\Delta ^\prime$.
\begin{remark}
  The ideas to prove Theorem \ref{th:generalized_BM_3.1} are
similar to those in \cite{bm}, with a few differences that improve
the probability of coalescence of 1/2 and 2/1 discrepancies. First,
  the labeling procedure in \cite{bm} was different: there,
{\rm all} $\eta^1$ particles were labelled (but none of the $\eta^2$
particles); they were called ``uncoupled'' when corresponding to 1/2
discrepancies, and ``coupled'' otherwise. Thus a coalescence of
discrepancies was called a ``coupling of labels''.   Secondly, we
introduce a notion of ``windows'' through stopping times slightly
different from those in \cite{bm}.
\end{remark}
 We want the processes $(X^{x,i}_t: t \geq 0)$ to be such  that \\
\indent
1)  for all $x \in (-LN,LN), i\in\{1,\cdots,K\}$,  if there are $h$  1/2 discrepancies at
$x$ at time $0$, then $X^{x,i}_0 = x$ for $i \leq h$, otherwise
$X^{x,i}_0 = \Delta,$\\
\indent
2) if $s<t$ and $X^{x,i}_s = \Delta $, then $X^{x,i}_t = \Delta $,\\
\indent
3) if there are $h$ 1/2 discrepancies at time $t$ at site $z$, then
there  exists  precisely $h$ pairs  $(x_j, i_j)$ so that $X^{x_j,i_j}_t = z$,  \\
\indent 4) for all  $(x,i)$ and $t$,  % relevant,
   the (random)  space-time point $X^{x,i}_t$
is either the position of a 1/2 discrepancy at time $t$ or equal to $\Delta$ and\\
\indent 5) for all $x \in \IZ$,   $i\in\{1,\cdots,K\}$,   $X^
{x,i}_. $ cannot jump except (possibly) at $t \in \N^{z,y} $ for some $z, y \in \IZ$ (neither of which may equal   $X^{x,i}_{t^-}$).   Equally, we insist that if some $t \in \N^{z,y}$ for some $z,y$ entails no change in both processes (i.e. $ \eta^1_{t^-} = \eta^1_t , \eta^2_{t^-} = \eta^2_t $), then there  will  be no movement of any of the $X^{x,i}_.$ processes at $t$.\\

Of course  for those five conditions to hold there can be many choices of the processes $(X^{x,i}_t: t \geq 0)_{x \in (-LN,LN), i \leq K}$.   We will make a choice that is natural, tractable and serves to prove Theorem \ref{th:generalized_BM_3.1}.

The choice of motions for the $X^{x,i}_.$ is ``solved'' for $p(.)$ a kernel of finite range  (see \cite{bm}).   For a general $p(.)$ we must be able to deal with  jumps  between sites $x $ and $y$ separated by great distances.  Accordingly we distinguish between changes in the $X^{x,i}_.$ processes occuring at $t \in \N^{y,z} $ for $|z-y|$ large and those contained in a Poisson process corresponding to a close pair  of sites.
We fix now an $\epsilon >0$ but arbitrarily small.  Associated with this $\epsilon $ we will choose an integer
$m = m_ \epsilon $ which will be large enough to satisfy various (increasing) properties that we will specify as our argument progresses.
The rules for the evolution of the $X^{x,i}_. $ at a point $t \in \cup_{y,z} \N^{z,y} $ will differ according to whether
$t \in \N^{z,y}$ for $|z-y| \geq m_\epsilon$  (we call such jumps ``big jumps'')  or not.
 We note that having finite systems of particles ensures that the rate at which relevant points in $\cup_{z,y}\N^{z,y}$ occur is bounded by  $K(2LN+1)||b||_\infty $.
Thus the time for jumps in the processes forms a discrete set, having no
cluster points.  Between these times we specify, by 5) above, that
$X^{x,i}_. $ must be constant for all  $x \in \IZ, i\in\{1,\cdots,K\}$.\\

We must now detail the motions of the $X^{x,i}_. $ at times
$t \in \N^{z,y}$.
As noted in 5) if no particle motion results then there is no motion of the discrepancies.  Furthermore if there are no 1/2 discrepancies at sites $z$ and $y $ then again no motion of 1/2 discrepancies results. Equally if at this instant a particle for each process moves from $z$ to $y$, then there is no motion of discrepancies.  This leaves two types of
 big jumps   occuring at $t$ to consider: $t \in \N^{z,y}$ for
$|z-y| \geq m_\epsilon$,  with a 1/2 discrepancy located either on $z$ or on $y$ at ``time'' $t^-$. \\

[a] A $\eta^1 $ particle moves from $z$ to $y$ (but not a $\eta^2 $ particle).  If at time $t^-$ there were no 1/2 discrepancies at $z$ then necessarily by  assumption %s {\em (A3)},
 {\em (A4)} we would have 2/1 discrepancies at $y$, thus 1/2 discrepancies neither on $z$ nor on $y$, a case we have excluded here. Therefore
 there are 1/2 discrepancies at $z$ at time $t^-$; we pick one at random, uniformly among  pairs
$(x,i)$ so that $X^{x,i}_{t^-} = z$ and move this discrepancy
 (and its label).

[a1]  If there are no $2/1$ discrepancies at $y $ at time $t^-$, then
the discrepancy chosen  and its label  jump to $y$;

[a2] if there are   $2/1$ discrepancies at $y $ at time $t^-$
then the $X^{x,i}_.$ chosen jumps to $\Delta $ and a  2/1
discrepancy is picked at random  at $y$  and  its label  jumps to $\Delta ^ \prime $ (and so each one must remain in these states thereafter: those 1/2 and 2/1 discrepancies have coalesced).  \\

[b]
A $\eta^2 $ particle moves from $z$ to $y$ (but not a $\eta^1 $ particle).

[b1] If at time $t^-$ there  are 2/1 discrepancies at $z$ we pick one at random, uniformly among these
 and move this discrepancy to $\Delta'$.
Since at time $t^-$ there are no 1/2 discrepancies on $z$, there must be some on $y$;
then one of these discrepancies is chosen uniformly at random and  its label
moves to $\Delta $.

[b2] If at time $t^-$ there is no 2/1 discrepancy at $z$, then (cf. [a] above) necessarily
 by assumption
 {\em (A4)}  there must exist 1/2 discrepancies at $y$. Again we choose one of these discrepancies at random and move it  (and its label)
 to $z$. \\

The motion of the $X^{x,i}_.$s for
$ t \in \N^{z,y}$ for some $|z-y| < m_ \epsilon $ is more complicated but follows along the lines of the rules introduced in \cite{bm}.\\

  We adopt an {\sl ordering}   $ \prec $   {\sl of   labels of
discrepancies} $X^{x,i}_. $ so that the spatial positioning is respected but which also orders  labels of  discrepancies on the same site.  The ordering among  ``active'' (in a sense made precise below)  discrepancies can only be changed by a   big
jump of size at least $m _ \epsilon $ for a 1/2 discrepancy
  (thus the jumps described in [a], [b] above),  at which point the   label of the  jumping discrepancy is assigned the lowest order among   labels of  1/2 discrepancies currently at the new site   (this choice is consistent with  the upper bound for $\Delta^{x,i}_t$ obtained below equation \eqref{def:delta}
  with respect to  the motions described in [a], [b], as will be explained later on).

Here a difference with the preceding cases is that at a single time $t$ many (but always a bounded number) $X^{x,i}_.$s may move so that labelled 1/2 discrepancies  keep their relative order. \\

[c] If at time $t^-$ neither site $z$ nor $y $ is the location of a 1/2 discrepancy
 then there is no motion for any $X^{x,i}_.$ at time $t$.

[d] If at time $t^-$ exactly one of the sites $z,y$ is the location
of 1/2 discrepancies, while the other site is not the current
position of 2/1 discrepancies,  then we fix the labels at time $t$
according to the following two requirements (we take $[z,y]$ to
signify $[y,z] $ in the case where $z$ exceeds $y$):  first
$X^{x,i}_t = X^{x,i}_{t^-}$ for all   pairs  $(x,i)$ for which
$X^{x,i}_{t^-} $ is outside $[z,y]$, secondly the $X^{x,i}_t$s are
chosen for   $X^{x,i}_{t^-} \in [z,y]$   so as to preserve order
(as in  \cite[Section 3]{bm}):   $X^{x,i}_{t^-} \prec X^{x^ \prime,
i ^ \prime}_{t^-} \Rightarrow X^{x,i}_{t} \prec X^{x^ \prime,i^
\prime}_{t}$.

[e] If at $t^-$ one of the  sites  $z,y$ is the location of  1/2 discrepancies and the other of $2/1 $ discrepancies,  then
 we relabel as follows:

[e1]  First we randomly select a random interval, called a
``window'' (see below) among the ``active windows" that contain both
$z$ and $y$.  Let this window be denoted $[u,v]$.  Then among all
pairs $(x,i)$ with $X^{x,i}_{t^-}\in[u,v]$ we choose (again all
candidates being equally likely) one $(x,i)$ and $X^{x,i}_t $ is
specified  to be $\Delta$, for the other pairs $(x^ \prime ,i ^
\prime ) $ we specify the $X^{x^\prime, i ^ \prime}_t$s so that
$X^{x^ \prime, i ^ \prime}_.$s outside $[u,v]$ remain where they
were while the order of $X^{x^\prime, i ^ \prime}_.$s within $[u,v]$
(apart from $X^{x,i}_.$) is preserved. Notice that this may result
in many (but a bounded number   of)   motions of labels: If e.g. the
motion is a 2/1 discrepancy at $z$ moving back to a 1/2 discrepancy
at $y$, but  $x$, the location of a 1/2 discrepancy whose label is
being chosen to be sent to $\Delta$ is such that $x
> z$,  then   labels of 1/2 discrepancies in $[y,x]$ are
shifted rightward  (or stay   on the same site if it is the location
of many labels).

[e2]   It may well be that the points $z$ and $y$ do not belong to a
single active window,   in which case $|z-y|<M_0$ or $|z-y| >
M_0+m_\epsilon$ (according to the definitions of $M_0$ and of
windows given below).   In this case the 1/2 discrepancy relevant to
the pair $z$ and $y$ at time $t$ has its label assigned to $\Delta$
and all other 1/2 discrepancies have their position (and label)
unchanged.
(Thus we are back to the behaviour described in [a],[b]).\\

\begin{remark}\label{awkward}
The relabeling enables us to get rid of the possibility of
a 2/1 discrepancy being close to a 1/2 discrepancy but not having a
chance of coalescing with it. Indeed, thanks to this manoeuvre,
whenever a 2/1 discrepancy comes close to a 1/2 discrepancy then
there is a nontrivial chance the label of the 1/2 discrepancy will
be sent to $\Delta$, while if we would have simply said that the
directly affected discrepancy has its label which goes to $\Delta$,
there might exist joint configurations where a 2/1 discrepancy is
close to a 1/2 discrepancy but the chance of it coalescing with that
particular discrepancy is essentially zero.
\end{remark}

It remains to describe the random intervals we call ``windows".  We follow closely the slightly different definition given in \cite{bm}.\\

In the following result a process on an interval $I$ will be a
process on state space $\{0,\cdots,K\}^{I}$ which obeys the same
evolution rules as  before, given  the Poisson processes $\N^{z,y}$ (and the
uniform random variables  $U(z,y,t)$  associated to $t\in\N^{z,y}$)
for $z,y \ \in \ I$.  We first observe that
 since by assumption {\em (A1)} kernel $p(.)$ is irreducible,
then for $n$ large enough \be \label{3} \mbox{ greatest common
divisor } \{x:p_n(x) \ne 0\} \ = \ 1, \ee where the (typically
sub Markov) kernel $p_n$ satisfies $p_n(x) = p(x)\indicator_{\{|x|
\leq n\}}$. The kernel $p_n(x)$ is finite range and we have as in
\cite[Lemma 3.1]{bm},

\begin{lemma} \label{lemfinite}
Let $n $ be sufficiently large that \eqref{3} holds.
For all $ m$ sufficiently large and all Harris coupled pairs of
 processes on $[0,m]$ evolving according to kernel $p_n(.)$  and $b(.,.)$,  $\eta^1_. $ and $\eta^2_.$ with  initial configurations
$\eta^1_0,\eta^2_0$ satisfying
$$
\eta^1_0(0) > \eta^2_0(0), \quad  \eta^2_0(m) > \eta^1_0(m) ,
$$
there is a strictly positive chance $c_{m}$ that there is a coalescence for the joint processes in  time interval $[0,1]$, that is that
$$
\sum_{x \in [0,m]}|\eta^1_1(x)-\eta^2_1(x)| < \sum_{x \in [0,m]}|\eta^1_0(x)-\eta^2_0(x)|.
$$
\end{lemma}

This immediately yields
\begin{cor} \label{corfin}
There exists $M_0$ so that for all $M \geq M_0$ if for
Harris coupled processes $\eta^1_., \ \eta^2_. $, for $0 \leq x \leq
x+M_0 \leq y \leq M, \eta^1_0(x) > \eta^2_0(x), \quad  \eta^2_0(y) >
\eta^1_0(y) $, then there is a strictly positive constant
$C_M$ so that with probability at least $C_M$ during time interval
$[0,1]$ (uniformly over all relevant joint initial configurations)

(i) there is no $t \in \N^{u,v} $ for $ u \in [0,M], v \notin [0,M]$ or vice versa,

(ii)
$$
\sum_{z \in [0,M]}|\eta^1_1(z)-\eta^2_1(z)| < \sum_{z \in [0,M]}|\eta^1_0(z)-\eta^2_0(z)|.
$$
\end{cor}

\proof Let $M_0 $ be a sufficiently large $m$ in the sense
of Lemma \ref{lemfinite} and $n$ be sufficiently large in the sense
of \eqref{3}.  Let $A$ be the event that in time interval $[0,1] $
there are no $t \in \N^{u,v} $ with either $ u \in [0,M], v \notin
[0,M]$ or $ u \in [x,y], v \notin [x,y]$, or vice versa.  Then
\[
\Prob(A)\ge e^{-4||b||_\infty \mu_1}
 \]
 where recall
$\mu_1 \ = \ \sum_w |w|p(w) < \infty $.  Furthermore event $A $ is independent of event
\[B=(
\mbox{there is no } t \in [0,1] \cap \N^{u,v}\mbox{ with }u,v \in [x,y]\mbox{ and }|u-v| \geq n),
\]
which has probability
\[\Prob(B)\ge e^{-||b||_\infty (M+1)\sum_{|w| \geq n} p(w)} \geq e^{-||b||_\infty(M+1)}.\]
 Conditional on
$A \cap B$, an event of probability
\[\Prob(A \cap B)\ge e^{-||b||_\infty (4 \mu_1 +M+1)}\]
the joint processes
$((\eta^1_s, \eta^2_s): 0 \leq s \leq 1)$ restricted to interval $[x,y]$ are just spatial translations of finite processes on $[0, y-x]$.  The result now follows from Lemma \ref{lemfinite}. \hfill$\Box$
\\

We now fix an $M_0$  (increasing $m_\epsilon$ if necessary),  so that $M_0 < m_\epsilon /10 $ and $M_0 > 10n$ where $n$ is sufficiently large in the sense of Lemma \ref{lemfinite}.
According to Corollary \ref{corfin}, the choice of $M_0$
is such that if a 1/2 discrepancy and a 2/1 discrepancy are separated by at
least $M_0$ (and less than $M_0 + m_\epsilon$) then there is a definite chance that
there will be a coalescence. If the separation is less than $M_0$, then, in
principle, we can say nothing about coalescence probabilities.

We are therefore ready to define ``windows", which will be space intervals of length $m_\epsilon +M_0$, on which coalescence will be favored.     A window will be associated to a  label of 1/2 discrepancy $X^{x,i}_.$.  Given
  $X^{x,i}_.$  (with $X^{x,i}_0\not=\Delta $)  we define the following stopping times:
 $T^{x,i}_0 = 0$,
\be\label{sigma} \sigma_{x,i} = \inf \{t\geq 0: X^{x,i}_t = \Delta\},\ee
  and for $j \geq 0$  (with the convention $\inf \emptyset=+\infty$)
\beq\label{Six-Tix}
S^{x,i}_j &=& \inf\{t \in [T^{x,i}_j,\sigma_{x,i}) : \exists \mbox { a 2/1 discrepancy in } \cr
 &&\qquad [X^{x,i}_t+M_0, X^{x,i}_t + M_0+m_ \epsilon ] \},\cr
a_j^{x,i}&=& X^{x,i}_{S^{x,i}_j},\cr
T^{x,i}_{j+1} &=& (S^{x,i}_j+1) \wedge \cr
 &&\inf\{t \geq S^{x,i}_j: t \in \N^{u,v} \mbox{ for  } u \in [a_j^{x,i} , a_j^{x,i} + M_0+m_ \epsilon ],\cr
 &&\qquad v \in [a_j^{x,i} , a_j^{x,i}+ M_0+m_ \epsilon ]^c \mbox{ or vice versa }\}
  \wedge \cr
  && \inf\{t \geq S^{x,i}_j: \sum _{u \in [a_j^{x,i}, a_j^{x,i} + M_0+m_ \epsilon ] } |\eta^1_t(u)- \eta^2_t(u)| < \cr
 &&\qquad\sum _{u \in [a_j^{x,i} , a_j^{x,i} + M_0+m_\epsilon ] } |\eta^1_{S^{x,i}_j}(u)- \eta^2_{S^{x,i}_j}(u)| \}
  \wedge \cr
  &&\inf \{t \geq S^{x,i}_j: \exists u, v \in [a_j^{x,i} , a_j^{x,i} + M_0+m_\epsilon ],\, t \in \N^{u,v} ,\cr
 &&\qquad
   |u-v| \geq m_\epsilon  \}.
\eeq
%
% (here if, for instance, for some $j$, between times $ T^{x,i}_j
%$ and $ \inf\{t \geq T^{x,i}_j: \exists \mbox { a 2/1 discrepancy in
%} [X^{x,i}_t+M_0, X^{x,i}_t + M_0+m_ \epsilon ] \}$ the stopping
%time $\sigma _{x,i}$ occurs, then $ S^{x,i}_j $ and all subsequent
%$S^{x,i}_., T^{x,i}_. $ are taken to be infinity).
   Times
$T^{x,i}_{j+1}$ and $S^{x,i}_j$ are defined in such a way that one
can use Corollary \ref{corfin} to conclude that there is a positive
chance of coalescence between times $S^{x,i}_j$ and $T^{x,i}_{j+1}$:
If the first or third event defining $T^{x,i}_{j+1}$  does not
occur then Corollary \ref{corfin} can be applied.  Since these two
events occur with finite rates we can expect coalescence with
positive probability between times $S^{x,i}_j$ and $T^{x,i}_{j+1}$.

For some  $(x,i)$ and $j$ with $S^{x,i}_j $ finite,  an  $((x,i),j)$
{\sl space window}  is a space interval $[a_j^{x,i} , a_j^{x,i} +
M_0+m_ \epsilon ]$.  It is taken to be {\sl active} during the  time
interval $[S^{x,i}_j, T^{x,i}_{j+1}]$,  called an $((x,i),j)$ {\sl
time window}. Indeed, the presence of a 2/1 discrepancy in
$[a_j^{x,i}+ M_0 , a_j^{x,i} + M_0+m_ \epsilon ]$ should favor a
coalescence with a 1/2 discrepancy  (cf. Remark \ref{awkward}). We
remark that a space window is only relevant while it is active, that
a point $u$ at a time $t $ may belong to several distinct  space
windows but that
this number is bounded by $M_0+m_\epsilon + 1 $, the size of a  space  window.\\

 We now define the evolution of the labels  $(Y^y_t : t \geq 0)$
of 2/1 discrepancies  (which will be more natural and intuitive than the processes of labels for 1/2 discrepancies).  Once a process
$Y^y_.$ hits $\Delta ^\prime $ it must remain at this ``position" ever after.
  We stipulate that
the $Y^y_.$ be a cadlag process which jumps at time $t$ only if for some $z \in \IZ, \ t \in \N^{z, Y^y_{t^-}}$ or
$t \in \N^{ Y^y_{t^-},z}$.  Furthermore nothing happens if at this time $t$ both a $\eta^1 $ and a $\eta^2$
particle move.

[f]
If for $t \in \N^{z, Y^y_{t^-}}$ solely a $\eta^1$ particle moves from $z$ to $Y^y_{t^-}$, one of the 2/1 discrepancies at $Y^y_{t^-}$ is randomly selected; then  if there were 1/2 discrepancies at $z$ at time $t^-$,  its label  moves to $\Delta ^\prime $ (this is case [a2] above when $|z-Y^y_{t^-}|\ge m_\epsilon$);  if there are no 1/2 discrepancies at $z$ at time $t^-$, it moves to $z$  as well as its label.    If for $t \in \N^{z, Y^y_{t^-}}$ solely a $\eta^2$ particle moves from $z$ nothing can happen to $Y^y_.$.

[g] If for  $t \in \N^{ Y^y_{t^-},z}$ only a $\eta^1$ particle
moves, then there is no motion for $Y^y_.$; if the  motion involves
uniquely a $\eta^2$ particle, then one of the 2/1 discrepancies
currently at site $Y^y_{t^- } $ is moved.  If at $t^- $ there is a
1/2 discrepancy at site $z$ then the  label of the  2/1 discrepancy
chosen moves to $\Delta ^ \prime $ (this is case   [b1]   above when
$|z-Y^y_{t^-}|\ge m_\epsilon$), if not  the 2/1 discrepancy (and its
label)   move to $z$.

Notice therefore that there is no relabeling scheme to preserve order for the processes $(Y^y_t : t \geq 0 )$.\\

\noindent
 To deal with Theorem \ref{th:generalized_BM_3.1}, we now consider the quantity \eqref{Phi}.
 Since by \eqref{1} the total number of particles is finite, $ \sup_{x\in\IZ} \Phi_t(x) $
 is equal to the maximum of $0$ and the  maximum over   $x\in\IZ,i\in\{1,\cdots,K\}$    of $\Phi_t(X^{x,i}_t)$.

We define for $(x,i),t$ such that $X^{x,i}_t\not=\Delta$,
 so that $\sigma_{x,i}>t$,
\beq\label{def:delta}\Delta_t^{x,i} &=&
%\sum_{y\geq X_t^{x,i}}(\eta^1_{t}(y) - \eta_{t}^2 (y)) - \sum_{y\geq X_0^{x,i}}(\eta_0^1(y) - \eta_0^2 (y)) \
 \widetilde \Phi_t(X^{x,i}_t) - \widetilde\Phi_0(X^{x,i}_0)\qquad\hbox{ with }\cr
\widetilde\Phi_t(X^{x,i}_t) &=&
\sum_{(y,j)\in\IZ\times\{1,\cdots,K\}}{\bf 1}_{\{X_t^{x,i}\prec X_t^{y,j}\}}\cr&&-\sum_{y\in\IZ: y> X_t^{x,i}}(\eta^1_{t}(y) - \eta_{t}^2 (y))^- .
\eeq
The quantity $ \widetilde\Phi_t(X^{x,i}_t) $ counts (the number of 1/2
discrepancies to the right of $X^{x,i}_t $ at time $t$ minus the
number of such 2/1 discrepancies) minus (the same quantity at time
0).  Thus $ \Delta_t^{x,i}$ is equal to the number of  labels of 1/2
discrepancies $X^{u,k}_. $ for which   $X^{u,k}_0 \prec X^{x,i}_0 $
but $ X^{x,i}_s \prec X^{u,k}_s$   for some $s \leq t$ and up to
time $t$ (that is, the labels of  1/2 discrepancies that appear in $
\widetilde\Phi_t(X^{x,i}_t) $ but were not in $\widetilde\Phi_0(X^{x,i}_0)$),  plus the
number of  labels of 2/1 discrepancies $Y^y_. $  that were ``in"
$\widetilde\Phi_0(X^{x,i}_0)$ but ``disappear" from $ \widetilde\Phi_t(X^{x,i}_t) $,
that is, those for which $Y^y_0 > X^{x,i}_0 $ but $Y^y_s < X^{x,i}_s
$ for some  $s \leq t$ (and up to time $t$) plus the number of
$Y^y_.$  (with $Y^y_0 > X^{x,i}_0 $) which jumped to $\Delta ^\prime
$ at a time $s \leq t$ so that  $X_{s^-}^{x,i} < Y^y_{s^-}$ and the
label $X^{u,k}_.$ of the 1/2 discrepancy which jumps to $\Delta $ at
time $s$ is such that $X^{u,k}_{s^-} \prec X^{x,i}_{s^-}$ (notice
that since $X^{x,i}_t\not=\Delta$, we have $(u,k)\not=(x,i)$).
%either
%\indent
%(i) there exists $(u,k)$ so that $X^{u,k}_{s^-} \prec X^{x,i}_{s^-}$ and $X^{u,k}_. $ jumps to $\Delta $ at time $s$, or
%
%\indent
%(ii) there exists $(u,k)$ so that $ X^{x,i}_{s^-}\prec X^{u,k}_{s^-} $, and
% at time $s$, $X^{u,k}_. $ jumps to $\Delta $ and $X^{x,i}_. $ is shifted rightward.
%
 \\

  Indeed, the ordering of labels of  1/2 discrepancies  for long jumps introduced earlier ensures that the upper bound for $ \Delta_t^{x,i}$ described above remains valid during such jumps.
 \\

 We now consider three classes of discrepancies contributing to
the above bound.  The first and second classes are not exclusive but
this does not concern us as we are interested in an upper bound for
$\Delta_t^{x,i}$.

\begin{description}%[(a)]
\item[(a)] $ \Delta_t^{x,i}(a) $ counts the number of 2/1 discrepancies that jump from $(X_s^{x,i} , \infty)$ to $ (-\infty, X_s^{x,i}) \cup \{ \Delta ^\prime \}$ for $s \leq t$ in some $\N^{u,v} $ with $|u-v| \geq m_\epsilon $ plus the number of  labels of  1/2 discrepancies $X^{w,k}_.$ for which there exists $s \leq t$ so that   $X^{w,k}_{s^-} \prec X^{x,i}_{s^-} =   X^{x,i}_s \prec X^{w,k}_s$.   It should be noted that necessarily, given the relabeling scheme in force, such a crossing must result from
a jump
of size greater than or equal to $m_\epsilon$ from $ (-\infty, X_s^{x,i})$ to $[X_s^{x,i} , \infty)$ by $X_.^{w,k}$.  Thus
$ \Delta_t^{x,i}(a) $ has only to do with big jumps.

\item[(b)] $\Delta_t^{x,i} (b) $ counts the number of  labels of  2/1 discrepancies, $Y^y_.$, for which for some  $s\in \N^{X^{x,i}_{s^-},u},s \leq t$,  for $u \geq X^{x,i}_{s^-} + m_\epsilon$ we have $X_{s^-}^{x,i} < Y^y_{s^-} = Y^y_s<  X_{s}^{x,i} $ plus the number of  labels of  1/2 discrepancies $X^{u,k}_.$ so that for some
 $s\leq t, \ X_s^{x,i} \prec X^{u,k}_{s^-} = X^{u,k}_s\prec  X_{s^-}^{x,i} $.   Again in the second case, given the relabeling scheme, it must hold at such an $s$ that  $|X_t^{x,i} - X_{t^-}^{x,i}|\geq m_\epsilon$
 %$X_s^{x,i} \leq X^{x,i}_{s^-} - m_\epsilon$
 (note in this case  $\partial\Delta_t^{x,i}:=\Delta_t^{x,i}-\Delta_{t^-}^{x,i} \leq  K \mid X_t^{x,i} - X_{t^-}^{x,i}\mid$).  Again,
$ \Delta_t^{x,i}(b) $ has only to do with big jumps.

\item[(c)] $ \Delta_t^{x,i}(c) $   deals with times $s \leq t$,   $s\in\N^{u,v}$ with $|u-v|< m_ \epsilon$, for motions described in case [e] above, with (in [e1]) or without (in [e2]) relabeling of some   1/2 discrepancies. The quantity  $ \Delta_t^{x,i}(c) $ counts the number of  labels of  2/1 discrepancies $Y^y_.$ which have not contributed to the two preceding random variables, so that   $Y^y_0 > X^{x,i}_0 $ and
at some time $s \leq t$,   $s\in\N^{u,v}$ with $|u-v|< m_ \epsilon$, $X_{s^-}^{x,i} < Y^y_{s^-}$ and either

\indent
(i) there exists  $(w,k)$  so that $X^{w,k}_{s^-}\prec X^{x,i}_{s^-}$ and at time $s$, $X^{w,k}_. $ jumps to $\Delta $ and $Y^y_s = \Delta ^\prime $
(this case includes both [e1] and [e2]), or\\
\indent
(ii) $X^{x,i}_. $  or $Y^y_.$  jumps at time $s$, and $Y^y_s<  X_{s}^{x,i} $
(if $X^{x,i}_. $ jumps, we are in case [d] above, hence there cannot be any 1/2 discrepancy between $X^{x,i}_{s^-} $  and $Y^y_{s^-}$).\\
(iii) there exists $(w,k)$ so that $ X^{x,i}_{s^-}<
Y^y_{s^-} < X^{w,k}_{s^-} $, and
 at time $s$, $X^{w,k}_. $ jumps to $\Delta $ and $X^{x,i}_. $ is shifted
 rightward.
 % $ \Delta_t^{x,i}(c) $ does not change. \\
\end{description}

So we have
$$
\Delta_t^{x,i}  \leq  \Delta_t^{x,i} (a) \ + \  \Delta_t^{x,i} (b) \  + \  \Delta_t^{x,i} (c) .
$$

We treat each term separately
\begin{lemma} \label{lemA}
There exists $c( \epsilon ) > 0$ so that for all $t$ sufficiently large
$$
\Prob(\Delta^{x,i}_t(a) > \epsilon t /5;\, \sigma_{x,i}>t) \ \leq \ e^{-c( \epsilon) t}.
$$
\end{lemma}

\proof
We can and will suppose that $m_ \epsilon $ has been fixed sufficiently large to ensure that
\be\label{vlarge}
\sum_{w\geq m_\epsilon }\hspace{0.2cm} w(p(w)+ p( -w))< \frac{\epsilon}{20||b||_\infty }.
\ee
The rate at which there is a jump of a 2/1 discrepancy from $(X_s^{x,i} , \infty)$ to $ (-\infty, X_s^{x,i}) \cup \{ \Delta ^\prime \}$ for $s \leq t$ in some  $\N^{u,y} $ with $|u-y| \geq m_\epsilon $
 is bounded by (because either solely a $\eta^2$ particle jumps, or
 solely a $\eta^1$ particle jumps that makes the 2/1 discrepancy move)
\beq
||b||_\infty \displaystyle\sum_{y> X_s^{x,i}} &\displaystyle{\sum_{u< X_s^{x,i}, u\le y-m_\epsilon} (p(y-u)+p(u-y))}\cr
&\displaystyle{\leq ||b||_\infty\sum_{w\geq m_\epsilon } w(p(w)+ p( -w))< \frac{\epsilon}{20}}
\eeq
   and similarly  for the rate for appropriate jumps of 1/2 discrepancies.
Thus these jumps are stochastically bounded by a rate $\epsilon /10 $ Poisson process.
So $\Prob(\Delta_{t}^{x,i} (a) > \epsilon t/5; \,\sigma_{x,i}>t) \leq e^{-c( \epsilon )t}$  for some $c > 0$ not depending on $N$.\hfill  $\Box$

\begin{lemma} \label{lemB}
There exists $c = c( \epsilon) > 0$ so that for all $N$ sufficiently
large and   $t \in [0,N]$
\[\Prob\left(|\{(x,i):   \Delta ^{x,i}_t(b)   \geq \frac{\epsilon N}{10} \} | \geq \frac{\epsilon^2}{50}N;\, \sigma_{x,i}>N \right) < e^{-cN}\]
 where $|A|$ denotes the cardinality of set $A$.
\end{lemma}

\proof    Since $\Delta ^{x,i}_t(b)$ is increasing in $t$ it is
sufficient to obtain the bound for $ t = N$.
 We use, for the moment, the fact that
\be\forall t \hspace{0.3cm} \partial \Delta_t^{x,i} (b) \leq K\vert X_{t}^{x,i} - X_{t^-}^{x,i} \vert\indicator_{\{\vert X_t^{x,i} - X_{t^-}^{x,i}\vert \geq m_\epsilon\}}
\ee
  though it should be noted that this gives a poor bound if the configurations $\eta_{t}^1$ and $\eta_{t}^2$ are ``close".
\vspace{0.2cm}

Observe that for  $(x,i)$ and  $(u,k)$  distinct, jumps of size
larger than  $m_\epsilon $ for $X^{x,i}_.$ and jumps for $X^{u,k}_.$
can be derived from independent Poisson processes of random but
bounded rates. Since the discrepancies are chosen uniformly
randomly when a Poisson clock rings at a site this claim is true for
two distinct discrepancies at the same site. Thus we can bound
stochastically the number of $(x,i)$ such that (recall
\eqref{vlarge})
\[ \displaystyle \sum_ {s\leq N}\vert X_s^{x,i} - X_{s^-}^{x,i} \vert \indicator_{\{\vert X_s^{x,i} - X_{s^-}^{x,i}
\vert \geq m_\epsilon\}} \geq \frac{\epsilon N}{10K}
\]
  by the number of $Z^{x,i}_.$ with $Z_N^{x,i} \geq (\epsilon
N)/(10K)$ where $Z_.^{x,i}$ are i.i.d. \hspace{0.15cm} random walks
all starting at zero which jump only in the positive direction
 by  $w \geq m_\epsilon$ at rate \hspace{0.15cm} $||b||_\infty (p(w) + p(-w))$.\\

If we chose $m_\epsilon$ sufficiently large then  for all  $ (x,i)$
\[\Prob\left(Z_N^{x,i} \geq \frac{\epsilon N}{10K}\right) \leq \frac{\epsilon^2}{100(2L +1)K}\]
for $N$ large by the law of large numbers and so we obtain
 \begin{eqnarray*}
&\displaystyle\Prob\left(|\{(x,i): \Delta_N^{x,i} (b) \geq \frac{\epsilon N}{10}\} |\geq \frac{N{\epsilon^2}}{50};\,\sigma_{x,i}>N\right)\\
\le & \displaystyle\Prob\left(\mbox{ Binom }((2L+1)N K,
\epsilon^2/100(2L+1)K) \geq \frac{N{\epsilon^2}}{50}\right) \leq
Ce^{-cN}.
\end{eqnarray*}
\hfill$\Box$

\vspace{0.2cm}

It remains to treat $\Delta_t^{x,i}(c)$.  As we have just seen, the variables $ \Delta_t^{x,i}(a)$ and $ \Delta_t^{x,i}(b)$
can be controlled by laws of large numbers applied to  big  jumps.  $\Delta_t^{x,i}(c)$ however is associated with the jumps of reasonable magnitude.
The 2/1 discrepancies contributing to $\Delta_t^{x,i} (c)$ should be split in two.  For a given pair $(v,k)$, we say a  time window
$[S^{v,k}_j, T^{v,k}_{j+1}]$  (or a $((v,k),j) $  space window)  is ``relevant'' to $X^{x,i}_.$ if the spatial interval  $[a_j^{v,k}, a_j^{v,k}+m_\epsilon +M_0]$ contains
$X^{x,i}_{S^{v,k}_j}$.
 We say $Y^y_.$ is ``associated" to a $((v,k),j) $  space window   $[a_j^{v,k}, a_j^{v,k}+m_\epsilon +M_0]$  relevant to $X^{x,i}_.$ if either of the following are true:\\
\indent
 1)  $Y^y_{S^{v,k}_j} \in [a_j^{v,k},a_j^{v,k}+m_\epsilon+M_0 ]$;\\
\indent
 2)  $Y^y_{T^{v,k}_{j+1}} \in [a_j^{v,k}, a_j^{v,k}+m_\epsilon +M_0]$  (notice that  if 2) occurs but not 1) then the $((v,k),j)$  space  window $ [a_j^{v,k},a_j^{v,k}+m_\epsilon+M_0 ]$ is desactivated by the entry of $Y^y_.$);\\
\indent
 3) $Y^y_{T^{v,k}_{j+1}} \in [ X^{x,i}_{T^{v,k}_{j+1}}+M_0, X^{x,i}_{T^{v,k}_{j+1}}+ m _\epsilon ]$.\\

 \noindent
 We say $Y^y_.$ is ``associated" if it is associated to one or more  space  windows.
 Otherwise $Y^y_.$ is not associated.  The sum over  labels of  2/1 discrepancies which are associated and contribute to $\Delta _t^{x,i}(c) $ is written $\Delta_t^{x,i}(c,{\rm ass})$ the contribution of nonassociated is $\Delta_t^{x,i}(c,{\rm non})$.\\

We note that no particle can leave or enter
$[a_j^{v,k},a_j^{v,k}+m_\epsilon+M_0 ]$ during the time interval
$[S^{v,k}_j, T_{j+1}^{v,k})$. Therefore the  2/1 discrepancies which
can contribute to $\Delta _t^{x,i}(c) $ during the time interval
$[S^{v,k}_j, T_{j+1}^{v,k})$  belong to
$[a_j^{v,k},a_j^{v,k}+m_\epsilon+M_0 ]$ and are located at sites to
the right of $X^{x,i}_{S^{v,k}_j}$ which implies that $\Delta
_t^{x,i}(c) $ can increase at most by $(M_0+m_\epsilon)K$ during
this time interval. Those 2/1 discrepancies which  are in  $[
X^{x,i}_{{(T^{v,k}_{j+1})}^-}, X^{x,i}_{({T^{v,k}_{j+1})}^-}+ m
_\epsilon ]$ at time ${(T^{v,k}_{j+1})}^-$ may lead to an increase
of at most $m_\epsilon K$ at time $T^{v,k}_{j+1}$. If
$Y^y_{T^{v,k}_{j+1}}\in[ X^{x,i}_{T^{v,k}_{j+1}},
X^{x,i}_{T^{v,k}_{j+1}}+ M_0)$, then $Y^y_s\in [ X^{x,i}_s+ M_0,
X^{x,i}_s+ M_0+ m _\epsilon ]$ for some $s<T^{v,k}_{j+1}$,   unless
it entered $[ X^{x,i}_{T^{v,k}_{j+1}}, X^{x,i}_{T^{v,k}_{j+1}}+
M_0)$ by a long jump of either $Y^y_.$ or $X^{x,i}_.$. As we can see
from case  3) above, any 2/1 discrepancy which enters  $[
X^{x,i}_{T^{v,k}_{j+1}}, X^{x,i}_{T^{v,k}_{j+1}}+ M_0)$ by a long
jump will not be associated.

Let $[S_j,T_j]$ $j \in \mathbb{N}$ denote the active windows
relevant to $(x,i)$ for all possible $(v,k)$ with reordered opening
times ($S_j \leq S_{j+1}$ but not necessarily $T_j \leq T_{j+1}$).
We are interested in finding an increasing subsequence of active
windows $[S_{j_k},T_{j_k})$ which are disjoint since we want to use
the strong Markov property to claim the independence of coalescence
events in such intervals to obtain our probability estimate.

Let $j_1 = 1$. Define for all $k \geq 1$,
$$j_{k+1} = \inf\{\ell > j_k:\,S_\ell \geq T_{j_k}\}.$$
Now we observe that $j_{k+1}\leq j_k + 2(M_0 + m_\epsilon)+1$ since
while the window $[S_{j_k},T_{j_k})$ is relevant to
$X^{x,i}_.$, then $X^{x,i}_.$ can belong to at most $2(M_0 + m_\epsilon)$ other
active windows. During the time interval $[S_{j_k},T_{j_k})$,
$X^{x,i}_.$ remains in the corresponding space window. Therefore
during this time interval $\Delta^{x,i}_.(c) $ increases by at most
$K(M_0+m_\epsilon)$. At time $T_{j_k}$ and time $S_{j_k}$,
$\Delta^{x,i}_.(c) $ increases by at most $m_\epsilon K$. No 2/1
discrepancy enters $[ X^{x,i}_t, X^{x,i}_t+ M_0)$ for $t\in
[T_{j_k},S_{j_{k+1}})$, unless it entered by a long jump (and
therefore is not associated). Also, no 2/1 discrepancy enters $[
X^{x,i}_t+ M_0, X^{x,i}_t+ M_0+ m_\epsilon ]$ (by a short or long
jump) during the same time interval. Thus we can conclude that
during the time interval $[S_{j_k}, S_{j_{k+1}})$, the total number
of associated 2/1 discrepancies which contribute to $\Delta^{x,i}_.(c)
$ is bounded above by $K\left(M_0+3m_\epsilon \right)$.
\\

\noindent
Treating $\Delta_t^{x,i} (c,{\rm ass})$ follows naturally along the same lines as with
 \cite[Proposition 3.2]{bm}:  during each relevant time window interval for label $X^{x,i}_.$,
 there is a reasonable probability that
$X^{x,i}_.$ jumps to  $\Delta $ and during such an interval
$[S_{j_k},S_{j_{k+1}})$ the number of $Y^y_.$ which are associated
is bounded by  $K(3m_\epsilon +M_0)$. Therefore if the event
$(\Delta_t^{x,i} (c,{\rm ass}) \geq \gamma N  ;\, \sigma_{x,i}>t)$
occurs for some $t$, then $X^{x,i}_.$ goes through at least $[\gamma
N]/(K(3m_\epsilon +M_0))$ successive disjoint time windows
$[S_{j_k}, S_{j_{k+1}})$ without jumping to $\Delta$ in the time
interval $[0,t]$. Since the probability of jumping to $\Delta$
during a given time window is equal to some $c'>0$, and jumps in
successive time windows are independent, the event $(\Delta_t^{x,i}
(c,{\rm ass}) \geq \gamma N ;\, \sigma_{x,i}>t)$ has a  probability
bounded above by
\[(1-c')^{[\gamma N]/(K(3m_\epsilon +M_0))}.\]
Thus we have in place of \cite[Proposition 3.2]{bm}
\vspace{0.2cm}

\begin{lemma} \label{lemC}
There exists $c, C \in (0, \infty)$ so that for all $t \geq 0$
\be \Prob(\Delta_t^{x,i} (c,{\rm ass}) \geq \gamma N  ;\,
\sigma_{x,i}>t)  \leq Ce^{-c \gamma N}.\ee
 \end{lemma}

It remains to assess  $\Delta_t^{x,i}(c,{\rm non})$   for $t \in
[0,N]$. Any $Y^y_.$ particle which is initially in $(X_0^{x,i},
X_0^{x,i} + M_0)$ may increase  $\Delta^{x,i}_t(c) $ without making a
long jump, if this occurs before the first time $X^{x,i}_.$ enters an
active window. This increase is bounded above by $K M_0$. Suppose
that $Y^y_.$ makes a contribution to this random variable. Then, by
definition, $Y^y_.$ is never associated with a space window relevant
to $X^{x,i}_.$. Since $Y^y_.$ does not contribute to $\Delta
_t^{x,i}(a) $ or $\Delta _t^{x,i}(b) $, it cannot traverse
$X^{x,i}_.$ via a jump of size greater than $m_\epsilon $, be it of
an $\eta^1$ or of an $\eta^2$ particle.

Thus the first time $s$ that $Y^y_s \in \ [X^{x,i}_s, X^{x,i}_s+
m_\epsilon ]$ must be less than   $t$.   Given that $Y^y_.$ is not
associated with a  space window relevant to $X^{x,i}_.$, it must be
the case that in fact at this point $s, \ Y^y_s \ \in \ [X^{x,i}_s,
X^{x,i}_s+ M_0 ] $ and that $s$ is the moment of a jump of size at
least $m_\epsilon $ either by $X^{x,i}_.$ or by $Y^y_.$.

{}From this one sees that    for all $ t \in [0,N]$   the
contribution $\Delta_t^{x,i} (c,{\rm non})$ is stochastically
bounded by a random variable that is in distribution the sum of
$KM_0$ times a Poisson random variable of parameter
$||b||_\infty t \sum_{z\geq m_\epsilon } p(z)$ and an independent
Poisson random variable of parameter $M_0||b||_\infty  t \sum_{z\geq
m_\epsilon } p(z)$.

\vspace{0.2cm}
Thus  for all   $ t \in [0,N]$,
$\Prob(\Delta_t^{x,i}(c) \geq \frac{\epsilon }{4} N; \sigma_{x,i}>N)
\leq C e^{-cN}$
 for some $c$ depending on $ \epsilon$ but not on $\eta_0^1, \eta_0^2$ or $N$.
\vspace{0.2cm}

So we have (recall  $\sum_{x\in\IZ} (\eta_0^1 (x)+ \eta_0^2 (x)) \leq 2K(2LN+1) )$,  by Lemmas \ref{lemA}, \ref{lemB} and \ref{lemC}

\begin{thm} \label{thm1}
For all $t \in [0,N]$,
\[\Prob\left(|\{(x,i) \hspace{0.2cm}: \hspace{0.2cm}\Delta _t^{x,i} \geq \frac{\epsilon N}{2}\} |\geq \epsilon^2 N; \sigma_{x,i}>N\right) \leq \hspace{0.2cm} Ce^{-cN}\]
  for $C, c \in (0, \infty)$ not depending on $N$.
\end{thm}
This result is close to the announced Theorem \ref{th:generalized_BM_3.1}.  The difference being that the latter deals with the supremum over all pairs $(x,i) $ of $ \Delta_t^{x,i} $, whereas Theorem \ref{thm1} asserts that the number of $(x,i)$ for which $ \Delta_t^{x,i} $ is ``too large" is small.\\
\\

\proof
 (of Theorem \ref{th:generalized_BM_3.1})

We prove the result by contradiction. By the conditions on the
initial configurations we necessarily have that $\sup_{x \in \IZ}
\Phi_0(x)$ is greater than or equal to zero.  So, for some $t \in
[0,N]$, $\sup_{x \in \IZ} \Phi_t(x)$ to exceed $\sup_{x \in \IZ}
\Phi_0(x)$ by some $\epsilon N$, we must have that for some pair
$(x,i)$ for a 1/2 discrepancy of label $X_t^{x,i}\not=\Delta$,
 \be \Delta _t^{x,i} \geq N \epsilon + \sup_{x \in \IZ}
\Phi_0(x). \ee
We have by Theorem \ref{thm1} that outside probability $Ce^{-cN}$,
\[|\{(v,k) :  \Delta_t^{v,k}  \geq \epsilon N/2 \}|\leq \epsilon^2 N.\]
Suppose for some $(x,i)$  with $X_t^{x,i}\not=\Delta$,
\[ \sum_{y\geq X_t^{x,i}}(\eta_{t}^1(y) - \eta_{t}^2(y)) -  \sum_{y\geq X_0^{x,i}}(\eta_{0}^1(y) - \eta_{0}^2(y)) >
\epsilon N.\]
This must mean that for at least  $(\epsilon N)/K$  sites $y\in
[X_t^{x,i}, +\infty)$, $\eta_t^1(y) > \eta_N^2(y)$. Let these points
be enumerated, in order, as $X_t^{x,i} \leq y_1 < y_2 \cdots < y_R$.
\vspace{0.2cm}

By hypothesis at most  $\epsilon^2 N (< ( \epsilon N)/(2K)$  for
$\epsilon$ small) are positions (at time $t$) of labels $X_t^{v,k}$
 with $X_t^{v,k}\not=\Delta$, and
$ \Delta_t^{v,k}  \geq \epsilon N/2$.  Let $i_1 = \inf \{j : y_j$ is
the position of a 1/2 discrepancy  label with $ \Delta_t^{v,k} \leq
\epsilon N/2\}$, then we have
\[ \sum_{y \geq X_t^{x,i}}\hspace{0.2cm} (\eta_t^1(y) -
\eta_t^2(y)) =
 \sum_{y =X_t^{x,i}}^{{y_{i_1}}-1}(\eta_t^1(y) - \eta_t^2(y)) {+} \sum _{y \geq y_{i_1}} (\eta_t^1(y) - \eta_t^2(y)).\]
But  $i_1 \leq \epsilon^2 N +1 $  so $\displaystyle
\sum_{y=X_t^{x,i}}^{y_{i_1}-1} (\eta_t^1(y) - \eta_t^2(y)) \leq
\epsilon^2 NK$, while
\[ \sum_{y \geq y_{i_1}}(\eta_t^1(y) - \eta_t^2(y)) \leq \frac{
\epsilon N}{2} + \sup_{w\in\IZ} \Phi_0(w) .\]
Thus $\displaystyle \sum_{y \geq X_t^{x,i}}(\eta_t^1(y) -
\eta_t^2(y)) < \epsilon N + \sup_{w\in\IZ}\Phi_0(w)$, a
contradiction. \hfill  $\Box$ \vspace{0.2cm}

 The above argument yields
\begin{cor}\label{cor2}
For $\eta_0^1$ and $\eta_0^2$ vacant on  $[-LN, LN]^c$
and satisfying for all $x \in [-LN, LN]$,
\[\vert \displaystyle \sum_{y\geq x} (\eta_0^1 (y) - \eta_0^2 (y)) \vert \leq \frac{\epsilon N}{4}, \]
we have outside probability $2Ce^{-cN}$ (for $C, c$ as in Theorem \ref{thm1})
\be \forall  t \in [0,N],\,\forall   x\in\IZ, \quad  \vert  \sum_{y\geq
x} (\eta_t^1 (y) - \eta_t^2 (y)) \vert \leq \epsilon N.\ee
\end{cor}
{}From our conclusions on finite configurations, one can compare
 infinite configurations through the following result for two initial
 configurations close on a (large) finite interval, one being finite
 (vacant outside that interval) and the other one infinite.
We wish to show that

\begin{thm}\label{thm2}  For $\epsilon >0$, there exists $L_0$ so that
for $L \geq L_0$, if $\eta_0^1$ is vacant on $[-LN, LN]^c$ and
$\eta_0^2$ satisfies
$$  \forall x \in [-LN, LN], \quad
 \vert \displaystyle \sum_{y=x}^{LN} (\eta_0^1 (y) - \eta_0^2 (y)) \vert \leq \frac{\epsilon N}{4},$$
 then for $C^1 , c^1$  depending on $L$, but not $N$, outside  probability $C^1 e^{-c^{1}N}$
 for all interval $I \subseteq  [-N,N]$ and for all $t \in [0,N]$,
$$ \vert \displaystyle \sum_{y \in I} (\eta_t^1(y) - \eta_t^2 (y))\vert \leq 3 \epsilon  N.$$
\end{thm}

\proof
We split
the process $\eta_.^2$ in  two classes of particles: first class particles  $\eta_.^{2.1}$  which at time $0$
were in the interval $[-LN, LN]$ and second class particles which are in $[-LN, LN]^c$ initially.
From Corollary  \ref{cor2}  we have that outside probability $C e^{-cN}$
\be \forall t \in [0,N],\, \forall x\in\IZ, \hspace{0.2cm} \vert
\displaystyle \sum_{y>x} (\eta _t^1 (y) - \eta_t^{2.1} (y)) \vert
\leq \epsilon N.\ee
This implies that for all $I \subset [-N,N]$, for all $t \in [0,N]$,
\be\vert \displaystyle \sum_{y \in I} (\eta_t^1 (y) - \eta_t^{2.1}
(y)) \vert \leq 2 \epsilon N.\ee
Thus to prove Theorem \ref{thm2} it will be enough to effectively bound
$$\displaystyle \sum_{y \in I} \eta_{t}^{2.2} (y) \leq \displaystyle \sum_{y \in [-N, N]} \eta_t^{2.2} (y),$$
the number of second class particle in $[-N,N]$ at time $t$.
\vspace{0.2cm}

The second class particles in $ [-N,N]$ can be divided in 2: those which jumped into $[-LN, LN]$ at the same
time as hitting $[-(L-1)N, (L-1)N]$, and those that enter $[-LN, LN]$ at a point in
$[-LN, -(L-1)N) \cup ((L-1)N, LN]$.
\vspace{0.2cm}

For the first we note that the entry of particles to $[-(L-1)N, (L-1)N]$ from $[-LN, LN]^c$  has (random) rate
bounded by

\be 2 \displaystyle \sum_{w\geq N}^{\infty} \sum_{y\geq w} ||b||_\infty(p(y) + p(-y))
\leq 2 \displaystyle ||b||_\infty\sum_{\vert w \vert \geq N} |w|p(w) \leq \frac{\epsilon}{10}\ee
 for $N$ large.
\vspace{0.2cm}

So the probability that the number of such entries over the time
interval $[0, t]$ exceeds $\epsilon /5 N$ is less than $H e^{-hN}$
for $H, h$ in $(0, \infty)$ not depending on $N$. \vspace{0.2cm}

For the remainder we note that the rate of the entrants to  $[-LN,
LN]$ must be bounded by $2\mu_1$ and so outside probability $H_{1}
e^{-h_1N}$ at most $ 4\mu_1N$  particles enter during time interval
$[0,t]$. \vspace{0.2cm}

If a second class particle enters  during time interval  $[0, t]$ at
$[-LN, LN] \setminus [-(L-1)N, (L-1)N]$   then for it to be in $[-N,
N]$ at time $t$, the sum of its absolute displacements over this
time interval must exceed  $(L-2)N$,  but for each such particle, %the sum of 
the absolute values of the jumps are stochastically
bounded by independent r.w.s.  $Z_N $  which jump over in positive
direction and jump to $w> 0$ at rate \hspace{0.2cm}
$||b||_\infty(p(w) + p( -w)) .$ \vspace{0.2cm}

Thus
\begin{eqnarray*}
&\Prob(| \{ \mbox{second class particles in }[-N, N]\mbox{ at time }t\mbox{ which entered }\\
&\mbox{ via }[-LN, LN] \setminus [-(L-1)N, (L-1)N]|) \} | \geq \frac{ \epsilon N}{4})\\
\leq &\Prob(\mbox{ Binom }(4\mu_1 N, p) > \epsilon  N/4) + H_1 e^{-h_1N}
\end{eqnarray*}
 where
\[p= \Prob(Z_N \geq (L-2)N) \rightarrow 0 \hspace{0.2cm} \mbox{ as }  N \rightarrow \infty\]
for  $L-2 > \mu_1$  (which is the case by \eqref{choice_L}).

\vspace{0.2cm}

Thus we obtain the defined bound for $N$ large, for all $t \in [0,N]$,
\[
\Prob\left(\displaystyle \sum_{y= -N}^{N} \eta_t^{2.2}(y) \geq
\epsilon N \right) \leq He^{-hN} + \varepsilon^{'} e^{-h1N }+
\varepsilon^{''} e^{-h''N} \leq Re^{-rN}
\]
and we are done. \hfill$\Box$
\section{Remaining lemmas}\label{macrostab2}
We need an extension of  \cite[Lemma 4.5]{rez} to nonfinite
range kernels:
\begin{lemma}\label{extended-Rezakhanlou}
Under the assumption $\mu_1=\sum_{z \in \IZ} |z| p(z) < \infty$, the measure $\nu^\rho$ has {\em a.s.} density $\rho$, that is
\be\label{eq:density_rho}
\lim_{l\to\infty}\frac{1}{2l+1}\sum_{x=-l}^l\eta(x)=\rho,
\quad \nu^\rho-\mbox{a.s.}
\ee
\end{lemma}

\proof
Let $\eta$ be a random
configuration with a distribution $\nu^\rho$. We want to show \eqref{eq:density_rho}.
We consider the stationary process $\eta_.$ with initial
distribution $\nu^\rho$. Since $\nu^\rho$ is translation invariant
we have
\[
\lim_{ l \to \infty} \sum_{x = -l}^l \frac{1}{2l+1} \eta_0(x) =
f(\eta_0)
\]
for $\nu^\rho$ almost every $\eta_0 \in {\bf X}$, where $f$ is a
translation invariant function. We will show that $f(\eta_0) =
f(\eta_1)$, $\nu^\rho \otimes \Prob$-a.s., thus showing that
$f$ is also a time invariant function. Let $n \in\IN$ and define the
event
\[
A^l_n = \{(\omega_0, \omega) \in (\Omega_0 \times \Omega): |\sum_{x
= -l}^l  \eta_0(\omega_0)(x) - \sum_{x = -l}^l
\eta_1(\eta_0(\omega_0),\omega)(x)|
> \frac{2l+1}{n}\}.
\]
 Now $|\sum_{x = -l}^l  \eta_0(\omega_0)(x) - \sum_{x = -l}^l
\eta_1(\eta_0(\omega_0),\omega)(x)|$ is the change in the net number of
particles in $[-l,l]$ during the time interval $[0,1]$. This can be written as
\[
(J^+_{-l}(\omega_0, \omega) + J^+_l(\omega_0, \omega)) -  (J^-_{-l}(\omega_0, \omega) + J^-_l(\omega_0, \omega))
\]
where $J^+_{-l}(\omega_0, \omega)$ is the total number of particles which jump from
$\IZ \cap (-\infty,-l-1]$ into $\IZ \cap [-l,l]$ during time interval $[0,1]$ and
$J^-_{-l}(\omega_0, \omega)$ is the total number of particles which jump from $\IZ \cap
[-l,l]$ into $\IZ \cap (-\infty,-l-1]$ during time interval $[0,1]$. We have that $J^\pm_l(\omega_0, \omega)$ are
defined similarly.  The sum  $(J^+_{-l}(\omega_0, \omega) + J^+_l(\omega_0, \omega))$ is bounded above by a
Poisson process with rate $||b||_\infty \mu_1$.
The same is true for $(J^-_{-l}(\omega_0, \omega) + J^-_l(\omega_0, \omega))$. Therefore $|\sum_{x =
-l}^l  \eta_0(\omega_0)(x) - \sum_{x = -l}^l
\eta_1(\eta_0(\omega_0),\omega)(x)|$ is bounded above by a Poisson random
variable with parameter $2 ||b||_\infty \mu_1$.
This implies that $\sum_{l=1}^\infty (\nu^\rho \otimes \Prob)
(A^l_n) < \infty$. Since this is true for all $n \in \IN$ we have
shown by Borel-Cantelli lemma that $f(\eta_0) = f(\eta_1),~ \nu^\rho \otimes \Prob$-a.s. Since $\nu^\rho \in {(\cal{I} \cap \cal{S})}_e$
this implies that $f$ is a constant, proving the result.
\hfill$\Box$
\\

  We finally state an extension to nonfinite range kernels of the {\em finite propagation} property at
particle level (see \cite[Lemma 5.2]{bgrs2}).
\begin{lemma}
\label{finite_prop-extended} There exist constant $v$,
and function $A(.)$ (satisfying $\sum_n A(n) < \ \infty $), depending only on $b(.,.)$ and
$p(.)$, such that the following holds. For any $x,y\in\IZ$, any
$(\eta_0,\xi_0)\in{\bf X}^2$, and any $0<t<(y-x)/(2v)$: if $\eta_0$
and $\xi_0$ coincide on the site interval $[x,y]$, then  with
$\Prob$-probability at least $1-A(t)$,
$\eta_s(\eta_0,\omega)$ and $\eta_s(\xi_0,\omega)$ coincide on the
site interval $[x+vt,y-vt]\cap\IZ$ for every $s\in[0,t]$.
\end{lemma}
\proof
Let  $\eta_0$ and $\xi_0$ be configurations in $\bf X$ which agree on all
sites $x$ such that  $m \leq x \leq n $. We couple the processes starting
from $\eta_0$ and $\xi_0$ by basic coupling so that they move together
whenever they can. We define random walks $L_t(\eta_t,\xi_t)$ and
$R_t(\eta_t,\xi_t)$  as follows.  Initially  $L_0  =m$ and $R_0 =n$. If $t\in\N^{z,w}$ for some sites $z < L_{t^-}\leq w
$, then $L_t = L_{t^-} + (w-L_{t^-} + 1)$. Similarly if  $t\in\N^{z',w'}$ for some sites $w' <L_{t^-}\leq
z'$, then $L_t = L_{t^-} + (z'-L_{t^-} + 1)$. Thus $L_t$ is a random
walk moving to the right. We define $R_t$ similarly moving to the
left. From the definition of $L_t$ and $R_t$ it follows that if, at
some time $t$, $a=L_t < R_t=b$, then during the time interval $[0,t]$
no particles entered or left $[a,b]$. Since  $\eta_0(x) = \xi_0(x)$
for all $a \leq x \leq b$ it follows that $\eta_t(x) = \xi_t(x)$ for
all $L_t \leq x \leq R_t$. Now the drift of $L_t$ can be written as
\begin{eqnarray*}
 v_L&=&||b||_\infty \left( \sum_{z < L_{t^-}} \sum_{w \geq L_{t^-}} (w-L_{t^-}+1) p(w-z)\right.\\
 \qquad&& + \left.\sum_{w' < L_{t^-}} \sum_{z' \geq L_{t^-}} (z'-L_{t^-} +1) p(w'-z')\right).
\end{eqnarray*}
 Using translation invariance of $p(.)$, then summation by parts,
the first term of $v_L$ can be written as
\[
||b||_\infty \sum_{j=1}^{+\infty}\sum_{i=j}^{+\infty}(i-j+1)p(i)=
||b||_\infty \sum_{i=1}^{+\infty} p(i) \sum_{j=1}^{i} j=
||b||_\infty \sum_{i=1}^{+\infty} \frac{i(i+1)}{2} p(i).
\]
  Similarly for the second term of   $v_L$   we write
\[
||b||_\infty \sum_{j=1}^{+\infty}\sum_{i=j}^{+\infty}jp(-i)=
||b||_\infty \sum_{i=1}^{+\infty} p(-i) \sum_{j=1}^{i} j=
||b||_\infty \sum_{i=1}^{+\infty} \frac{i(i+1)}{2} p(-i).
\]
Both terms are finite because  of the moment assumptions on $p(.)$, and
\[
v_L=||b||_\infty \frac{\mu_1+\mu_2}{2}
\]
where $\mu_2=\sum_{z\in\IZ} z^2p(z)$.
We can proceed similarly with $R_t(\eta_t,\xi_t)$  to show that the drift
$v_R$ of $R_t$ is $v_R = -v_L$. From the argument above it follows
that if
\[
\sum_{z\in\IZ} |z|^k p(z) < \infty
\]
for $k > 1$ then both $L_t$  and $R_t$ have finite   $(k-1)$th
moment. Since we have assumed that $p(.)$ has finite third  moment,
we can conclude that $R_t$ and $L_t$ have finite second moment.

Therefore if we take $v$ bigger than $v_L= -v_R$, then for all $t\ge
0$, $\Prob\{L_t\geq v_Lt\}\leq  A([t])$, for some function $A$
satisfying the announced finiteness condition which depends only on
$p(.)$, $b(.,.)$  and   $v- v_L$. \hfill$\Box$

\mbox{}\\ \\
\section*{Acknowledgements}   We thank the referee for a
careful reading and useful comments.   The research of T.S.M. is
partially supported by the SNSF, grants $\#\, 200021-107425/1$ and
$\#\, 200021-107475/1$. K.R. was supported by NSF grant DMS 0104278.
E.S. thanks EPFL for hospitality.  Part of this work was done during
the authors'stay at Institut Henri Poincar\'e, Centre Emile Borel
(whose hospitality is acknowledged), where T.S.M. and K.R. benefited
from a CNRS ``Poste Rouge'' for the semester ``Interacting Particle
Systems, Statistical Mechanics and Probability Theory''.

\end{document}